\documentclass[10pt,twoside]{article}
\usepackage{amssymb}
\usepackage{fancyhdr}
\usepackage{titlesec}
\usepackage{cite}
\usepackage{ifthen}
\usepackage{graphicx}

\titleformat{\section}{\centering\large\bfseries}{\S\arabic{section}}{1em}{}
\newboolean{first}
\setboolean{first}{true}

\begin{document}

\title{\bf \Large The Analysis of Rotated Vector Field for the Pendulum\author{Lian-Gang Li
 \\ \small \emph{School of Sciences, Tianjin University, Tianjin 300072, China}}\date{}}
\maketitle

\setlength\abovedisplayskip{2pt}
\setlength\abovedisplayshortskip{0pt}
\setlength\belowdisplayskip{2pt}
\setlength\belowdisplayshortskip{0pt}

\footnote{\emph{E-mail addresses}: liliank@tju.edu.cn}

\begin{center}
\begin{minipage}{135mm}
{\bf \small Abstract}.\hskip 2mm {\small The pendulum, in the presence of
linear dissipation and a constant torque, is a non-integrable, nonlinear differential
equation. In this paper, using the idea of rotated vector fields,
derives the relation between the applied force $\beta$ and the periodic
solution, and a conclusion that the critical value of $\beta$
is a fixed one in the over damping situation. These results are of practical significance
in the study of charge-density waves in physics.\vskip 5mm \bf \small
MSC}:\hskip 2mm {\small 34C25, 37C27, 34C05}
\end{minipage}
\end{center}

\vskip 6mm

\setcounter{section}{1} \setcounter{equation}{0}
\renewcommand{\theequation}
{\arabic{equation}}
 \setcounter{figure}{0}\renewcommand{\thefigure}
{\arabic{figure}}
\begin{flushleft}
\large \bf 1.\hspace*{2mm}Introduction\\
\end{flushleft}

The pendulum is a famous physical model due to its versatile roles in the history of
physics and mathematics\cite{Mawhin,Ort1,Ort2,Ahmad}, and not only interesting
for being the familiar example of nonlinear dynamics, but more important
for it manifesting in many important areas of physics, for example, charge-density wave (CDW),
a special phenomenon of electronics condensing in the low dimensional materials \cite{Fleming},
such as is also observed in the conventional 
 and the High Temperature Superconductor materials \cite{Morosan,Seib,Zabo}. To describe these special
properties in the electric conduction with CDW, Gr\"{u}ner and his co-authors 
proposed a single particle model \cite{Gru}, and come up with a
equation as following:
\begin{equation} \label{eq:f1}
\frac{\mathrm{d}^{2}\phi}{\mathrm{d}t^{2}}+\Gamma\frac{\mathrm{d}\phi}{\mathrm{d}t}+\sin\phi=\beta,
\end{equation}
where $\beta=\frac{E}{E_{0}}$, $E$ is an electrical Field applied,
$E_{0}$ is a constant, $\Gamma$ is the friction coefficient.
Equation (\ref{eq:f1}) is the mathematical model of a planar
pendulum in the presence of linear dissipation and a constant torque,
also called Gr\"{u}ner Equation in the CDW theory.
Due to being a non-integrable equation, Gr\"{u}ner just gave an approximate solution with neglect of
$\frac{d^{2}\phi}{dt^{2}}$ term \cite{Gru}. Although, an article gave the description for its
solutions with the analysis of vector fields rigorously [6],
it is merely an scenario, considering $\beta$ as fixed values, to describe the behavior
of its solutions with the $\Gamma$ changing. By the idea of \emph{rotated} vector fields,
it is the purpose of this paper, to elucidate the behavior of its solutions with the $\beta$ changing,
which just meets the needs for clarifying the mechanism of CDW in \cite{Li}.
\par
In 1881-1886, Poincar\'{e} established a theory named "vector field" in his excellent papers \cite{Poinca}.
With this approach, the solutions of differential equation can be regarded as integral
curves in the space of phase, and the qualitative properties of solutions
can be obtained with geometrical way. After that, G. F. Duff proposed the notion of the rotated vector
fields \cite{Duff}, then G. Seifert et al. developed it to that of the general
rotated vector fields \cite{Seif,Chen,Ma}, in which the essence of the idea, i.e. \emph{rotating}, is also employed in this paper.
\par
By the idea of rotated vector fields, this paper proves the
relation between $\beta$ and the periodic solution of the
equation, and find that if the $\Gamma$ is
bigger enough (over damping), the critical value of $\beta$, with
which this equation will have a critical periodic solution, will be a
fixed value all the time, namely, $\beta_{0}\equiv 1$.
These conclusions serve as essential roles in \cite{Li}.

\setcounter{section}{1} 
\renewcommand{\theequation}
{\arabic{equation}}
\begin{flushleft}
\large \bf 2.\hspace*{2mm}Preliminary definitions and lemmas\\
\end{flushleft}

\par
As traditionally done in dynamics analysis, after being written as a system of first-order equations, Equation (\ref{eq:f1}) is embedded into the $\phi-z$ phase plane $\mathbb{R}^2$:
\begin{equation}\label{eq:f2}
 \left\{
\begin{array}{ll}
\frac{\mathrm{d}\phi}{\mathrm{d}t}=z\\
\frac{\mathrm{d}z}{\mathrm{d}t}=\beta-\sin\phi-\Gamma\cdot z
\end{array}\right.,
\end{equation}
where $\Gamma>0$, $\beta\geqslant0$. In this situation a solution of
equation (\ref{eq:f2}) will correspond to a trajectory on the phase plane. According
to $\frac{d\phi}{dt}=z$, shows that the trajectories direct from
the left to the right while $z>0$, and they direct from the right to
the left while $z<0$.
\par
When $\beta>1$, there exists no equilibrium point. When
$0\leqslant\beta\leqslant1$, there exist equilibria on the $\Phi$-axis,
which coordinates are denoted by $(\phi_{n},0)$, where
\begin{equation}\label{eq:f3}
 \left\{
\begin{array}{ll}
\phi_{n}=n\pi+(-1)^{n}\phi_{0}\\
\phi_{0}=\arcsin\beta
\end{array}\right.,  \quad\quad n\in\mathbb{Z}.
\end{equation}
The equilibria have the properties as
following, and denoted by the following symbols as default.
 \vskip 2mm
\begin{minipage} [32mm]{138mm}
1)\hskip 2mm When $0\leqslant\beta<1$, $A_{k}$ at $(\phi_{2k-1},0)$, for $k\in \mathbb{Z}$, denote the
saddles of equation (\ref{eq:f2}). The slopes of two separatrices, passing through each of the saddles, are
respectively
\begin{equation} \label{eq:f4}
\lambda_{1}=-\frac{\Gamma}{2}+\sqrt{\left(\frac{\Gamma}{2}\right)^{2}+\cos\phi_{0}}
\hskip 1mm ,
\end{equation}
\begin{equation} \label{eq:f5}
\lambda_{2}=-\frac{\Gamma}{2}-\sqrt{\left(\frac{\Gamma}{2}\right)^{2}+\cos\phi_{0}}
\hskip 1mm .
\end{equation}
When $\Gamma>0$, $B_{k}$ at $(\phi_{2k},0)$, for $k\in \mathbb{Z}$, denote the stable focus.
\par
When $\Gamma=0$, $B_{k}$ at $(\phi_{2k},0)$, for $k\in \mathbb{Z}$, denote the centers. But we don't concern this case.
\par
2)\hskip 1mm When $\beta=1$, $\phi_{0}=\frac{\pi}{2}$,
$\phi_{2k+1}=\phi_{2k}$, namely, $A_{k+1}$ and $B_{k}$ coincide with
each other in position. They become combining equilibria.
 \par
\end{minipage}
\vskip 2mm 
Obviously,
equation (\ref{eq:f2}) possesses a period of $2\pi$ along the $\Phi$-axis,
so it can be treated as the cylinder $E^2=S^1\times\mathbb{R}^1$.
Due to the periodicity of $\phi$ on $E^2$, we just need to discuss the
interval $[-\pi-\phi_{0},\pi-\phi_{0}]$ of $\phi$ without special declaiming
and name it the main interval.
\par
\begin{figure} [h]
\centering
\includegraphics[width=2.80in, height=1.80in]{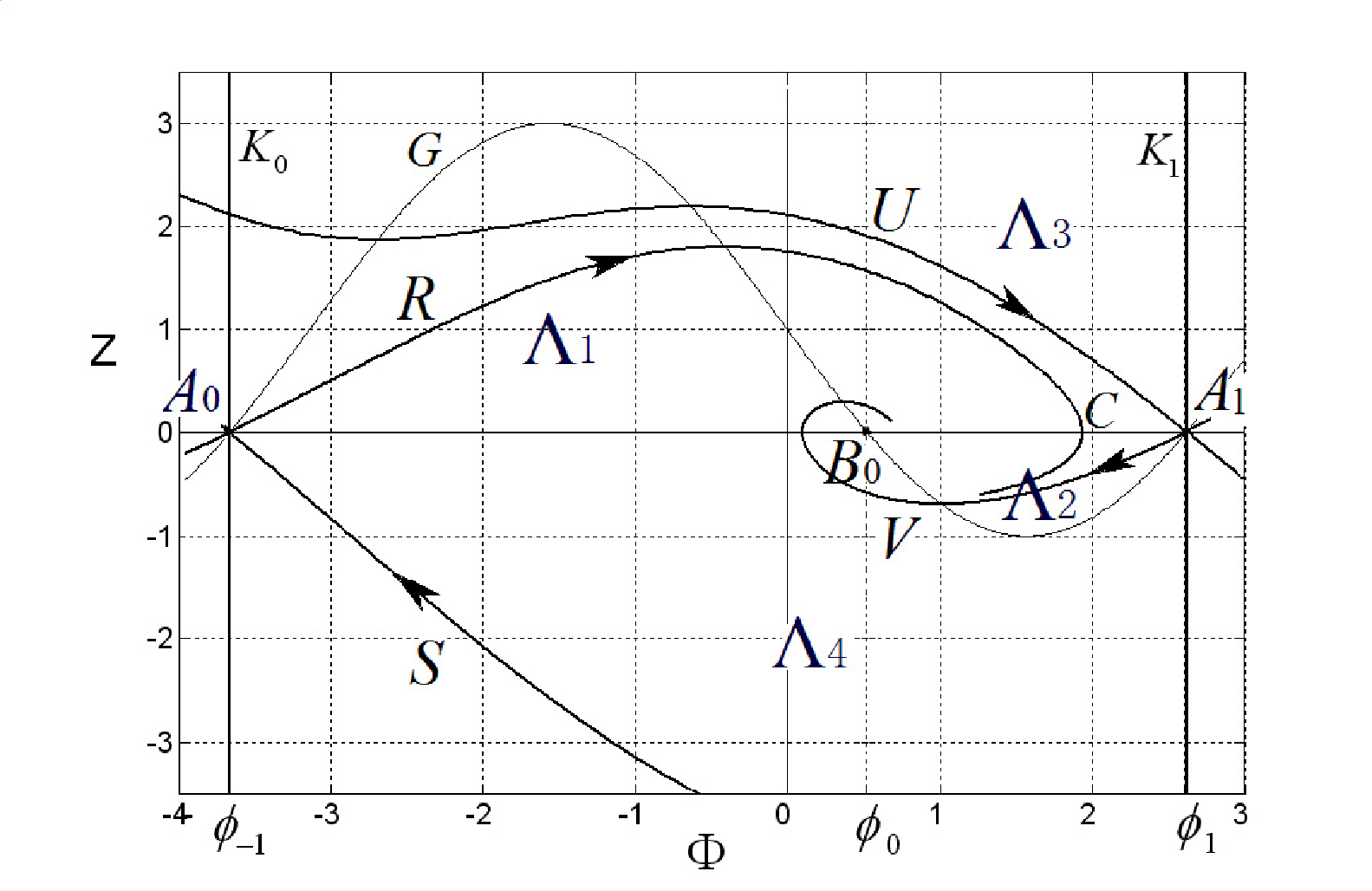}
\caption{Special trajectories in the main interval.}
\label{fig:p1}
\end{figure}

While $\Gamma>0$, $0\leqslant\beta<1$, in the main interval,
Let $R$, $V$ denote the trajectory emitting from $A_{0}$ and $A_{1}$
respectively; Let $S$, $U$ denote the trajectory coming into $A_{0}$
and $A_{1}$ respectively; as shown in Figure \ref{fig:p1}. In order to demarcate different regions,
we draw a curve $G$, a nullcline, as shown in Figure \ref{fig:p1}, which reads

\begin{equation} \label{eq:f6}
z=\frac{\beta-\sin\phi}{\Gamma}.
\end{equation}
\par\noindent
\par
\begin{minipage} [32mm]{138mm}
 Let $\Lambda_{1}$ denote the region encircled by $G$ and the $\Phi$-axis for
$z>0$ ;\\
 Let $\Lambda_{2}$ denote the region encircled by $G$ and the $\Phi$-axis for $z<0$;\\
 Let $\Lambda_{3}$ denote the region above $G$ and
$z>0$;\\
 Let $\Lambda_{4}$ denote the region below $G$ and $z<0$.\\
\end{minipage}
\par
 From equation (\ref{eq:f2}), it is easy to obtain:
\begin{equation} \label{eq:f7}
\frac{\mathrm{d}z}{\mathrm{d}\phi}=\frac{\beta-\sin\phi}{z}-\Gamma.
\end{equation}
\par
   Obviously, from equation (\ref{eq:f7}) it follows the properties as below:
\par
\begin{minipage} [32mm]{138mm}
a) $\frac{\mathrm{d}z}{d\phi}>0$, while $(\phi,z)\in\Lambda_{1}.or.\Lambda_{2}$;\\
b) $\frac{\mathrm{d}z}{d\phi}<0$, while $(\phi,z)\in\Lambda_{3}.or.\Lambda_{4}$;\\
c) $\frac{\mathrm{d}z}{d\phi}=0$, while $(\phi,z)\in G$;\\
d) $\frac{\mathrm{d}z}{d\phi}=\infty$, while z=0 and
$\phi\neq\phi_{n}$.
\end{minipage}
\par
 According to the above properties, it is exactly straightforward to prove the Lemma 1 as follows.

\par\noindent
{\bf Lemma 1.} {\it Let $K_{1}:\phi=\phi_{1}$ denote a line
perpendicular to the $\Phi$-axis through point $A_{1}$, as  shown in
Figure 1. Then the trajectory $R$ emitting from $A_{0}$ should passe
over the upper half of the $\Phi$-axis and must reach point $C$ at $(\phi_{c},0)$
of the $\Phi$-axis where $\phi_{0}\leqslant\phi_{c}\leqslant\phi_{1}$, or
at $(\phi_{1},z_{c})$ of line $K_{1}$ where
$0<z_{c}<\frac{\beta+1}{\Gamma}$.}

we can also give Lemma 2.
\par\noindent
{\bf Lemma 2.} {\it If equation (\ref{eq:f2}) exists a periodic
solution $z=Z(\phi)$, where $Z(\phi)=Z(\phi+2\pi)$, $-\infty<\phi<+\infty$, then it
must satisfy:
\begin{equation} \label{eq:f8}
\int_{\phi}^{\phi+2\pi}Z(\phi)\,\mathrm{d}\phi=\frac{2\pi\beta}{\Gamma}.
\end{equation} }
The proof of the Lemma 2 is very simple. Integrating
equality (\ref{eq:f7}) from $\phi$ to $\phi+2\pi$ and using its
periodicity, equality (\ref{eq:f8}) is obtained easily.
\par
For the proof in the next section, we also provide two definitions
as following.
\par\noindent
{\bf Definition 1.} {\it Here $\beta$ is fixed, take
$\Gamma$ as a parameter, $\theta$ is the angle between the vector and the
$\Phi$-axis, the changing ratio of $\theta$ with respect to $\Gamma$
is
 \begin{equation}\label{eq:f9}
 \frac{\mathrm{d}\theta}{\mathrm{d}\Gamma}=-\frac{z^{2}}{z^{2}+(\beta-\sin\phi-\Gamma \cdot
 z)^2}\leqslant0\quad .
 \end{equation}
 }
(It constructs a rotated vector field in terms of the notion in \cite{Duff}.)
\par\noindent
{\bf Definition 2.} {\it Here $\Gamma$ is fixed,
take $\beta$ as a parameter, $\theta$ is the angle between the vector and
the $\Phi$-axis, the changing ratio of $\theta$ with respect to
$\beta$ is
 \begin{equation}\label{eq:f10}
 \frac{\mathrm{d}\theta}{\mathrm{d}\beta}=\frac{z}{z^{2}+(\beta-\sin\phi-\Gamma \cdot
 z)^2}\quad .
 \end{equation}
}
(It constructs a general rotated vector field in terms of the notion in \cite{Seif,Chen,Ma}.)
\par\noindent
Notice, in the arguments, that of $z>0$  is involved only,
and $-\frac{\pi}{2}<\theta\leqslant\frac{\pi}{2}$ with $\beta$ or $\Gamma$ changing.
\par
Now, we give the third lemma.
\par\noindent {\bf Lemma 3.} {\it Here is $\beta=0$
, for any $\Gamma_{1}>0$, while $\Gamma=\Gamma_{1}$, then the
trajectory $R$ emitting from $A_{0}$ should passes over the upper half of
plane and must reach point $C$ at $(\phi_{c},0)$ of the
$\Phi$-axis, where $\phi_{0}\leqslant\phi_{c}<\phi_{1}$, i.e.,
between focus $B_{0}$ and saddle $A_{1}$(excluding $A_{1}$).}
\par
\begin{figure} [h]
\centering
\includegraphics[width=2.80in, height=1.80in]{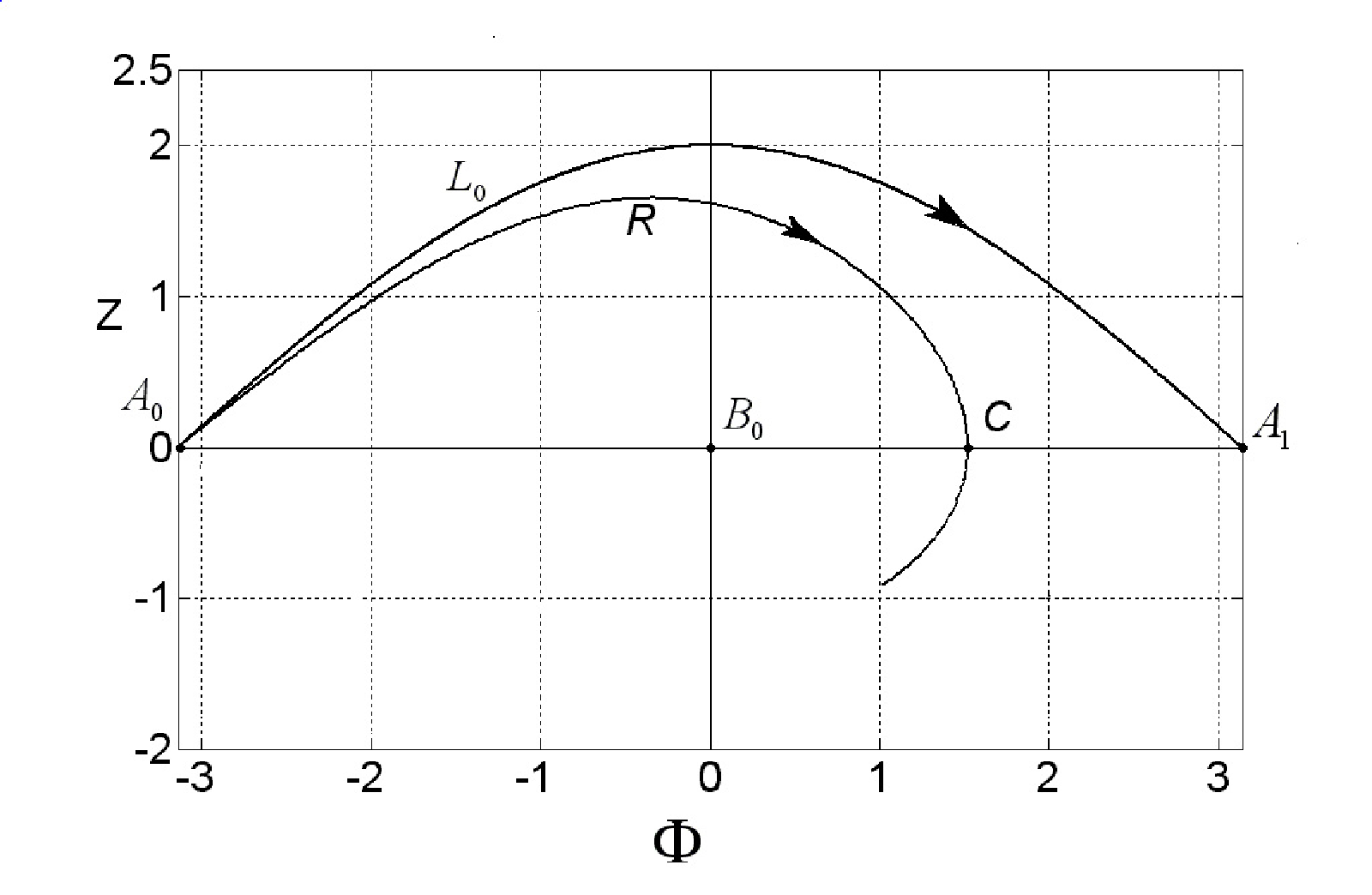}
\caption{Trajectories as $\beta=0$.}
 \label{fig:p2}
\end{figure}

\noindent
{\bf Proof.} Let $\Gamma=0$, $\beta=0$, thus equation
(\ref{eq:f2}) degenerates into a integrable system, sine-Gordon's Equation,
namely, ${\displaystyle\frac{\mathrm{d}^{2}\phi}{\mathrm{d}t^{2}}+\sin\phi=0}$.
It is easy to give a solution,
\begin{equation}\label{eq:f11}
z=\sqrt{2(\cos\phi+1)}.
\end{equation}
As shown in Figure \ref{fig:p2}, it corresponds a trajectory $L_{0}$
connecting the saddle $A_{0}$ with $A_{1}$.
\par\noindent
The slope of $L_0$ at saddles $A_{0}$ and $A_{1}$ is $\pm1$
respectively. Let $\Gamma =\Gamma_{1}$ and consider $\Gamma_{1}>0$,
$\beta=0$, from equality (\ref{eq:f4}) it follows that the
slope of trajectory $R$ at point $A_{0}$ is $\lambda_{1}$, where
$0<\lambda_{1}<1$. Therefore, $R$ is at the inside
of curve $L_{0}$ at the beginning. From Definition 1,
it follows that, except points $A_{0}$ and $A_{1}$,
$\frac{\mathrm{d}\theta}{\mathrm{d}\Gamma}<0$ for the vectors of field on curve $L_{0}$,
which means them rotated clockwise while $\Gamma =\Gamma_{1}$,
moreover $-\frac{\pi}{2}<\theta\leqslant\frac{\pi}{2}$, so the vectors
on curve $L_{0}$ point to the inside of $L_{0}$. From the
uniqueness of solutions, it follows that the trajectory $R$ runs at the inside of
curve $L_{0}$ as shown in Figure \ref{fig:p2}, or intersects at saddle $A_{1}$.
Assuming that trajectory $R$ intersects at the point $A_{1}$,
thus, the slope of $R$ at the point $A_{1}$ is $\geqslant-1$ as it
approaches the point $A_{1}$ from the left inside of $L_{0}$. However, it follows from equality (\ref{eq:f5})
its slope is $\lambda_{2}<-1$. This contradicts the assumption. So $R$ can not reach the point $A_{1}$.
According to Lemma 1, one can declare that $R$ must passes over the upper
half plane and reaches the point $C$ at $(\phi_{c},0)$ on the $\Phi$-axis,
where $\phi_{0}\leqslant\phi_{c}<\phi_{1}$, i.e., between point $B_{0}$ and
$A_{1}$(excluding $A_{1}$) .
\par The proof is complete.

\setcounter{section}{1} 
\renewcommand{\theequation}
{\arabic{equation}}
\begin{flushleft}
\large \bf 3.\hspace*{2mm} The proof of theorem\\
\end{flushleft}
\par
There is no cycle of type I on the cylinder $E^2$ for the system as Bendixson's Criteria \cite{Perk}, so in this paper,
the "periodic trajectory" or "periodic solution", refers to the cycle of
type II on $E^2$ as default, in some papers, also called the running periodic trajectory on $\mathbb{R}^2$.
\par
Now we propose the theorem as following, and use the preceding lemmas to prove it.
\par
\par \noindent
{\bf Theorem.} {\it Suppose $\Gamma =\Gamma_{1}$, consider $\Gamma_{1}>0$,
$R(\beta)$ denotes the trajectory emitting from point $A_{0}$, the following propositions hold.}
\par
\vskip 2mm {\it
\begin{minipage} [32mm]{140mm}
a)\hspace*{2mm} Equation (\ref{eq:f2}) must have a value of
$\beta_{0}(\Gamma_{1})$, $0<\beta_{0}\leqslant1$, while
$\beta=\beta_{0}(\Gamma_{1})$, $R(\beta_{0})$ must intersect
the $\Phi$-axis at point $A_{1}$. Here this trajectory is named as
critical periodic solution or pseudo-periodic solution. $\beta_{0}$ is called
the critical value, with which the equation will have a
critical periodic solution.
\par
b)\hspace*{2mm} $\beta_{0}(\Gamma_{1})$ is unique.
It means that there is only one value of $\beta_{0}$ for a given $\Gamma_{1}$.
\par
c)\hspace*{2mm} While $\beta>\beta_{0}$, equation (\ref{eq:f2})
exists a periodic solution $z=Z_{T}(\phi)$, that is
$Z_{T}(\phi)=Z_{T}(\phi+2\pi)$, and $Z_{T}(\phi)>0$, $-\infty<\phi<+\infty$;
while $\beta<\beta_{0}$, equation (\ref{eq:f2})
does not exist any periodic solution.
\par
d)\hspace*{2mm} The periodic solution of equation (\ref{eq:f2}) is
unique and stable.
\par
e)\hspace*{2mm} Suppose $\beta_{0}(\Gamma_{1})=1$, and if
$\Gamma_{2}>\Gamma_{1}$, then $\beta_{0}(\Gamma_{2})=1$.
\end{minipage} }
\par
\begin{figure} [h]
\centering
\includegraphics[width=2.80in, height=1.79in]{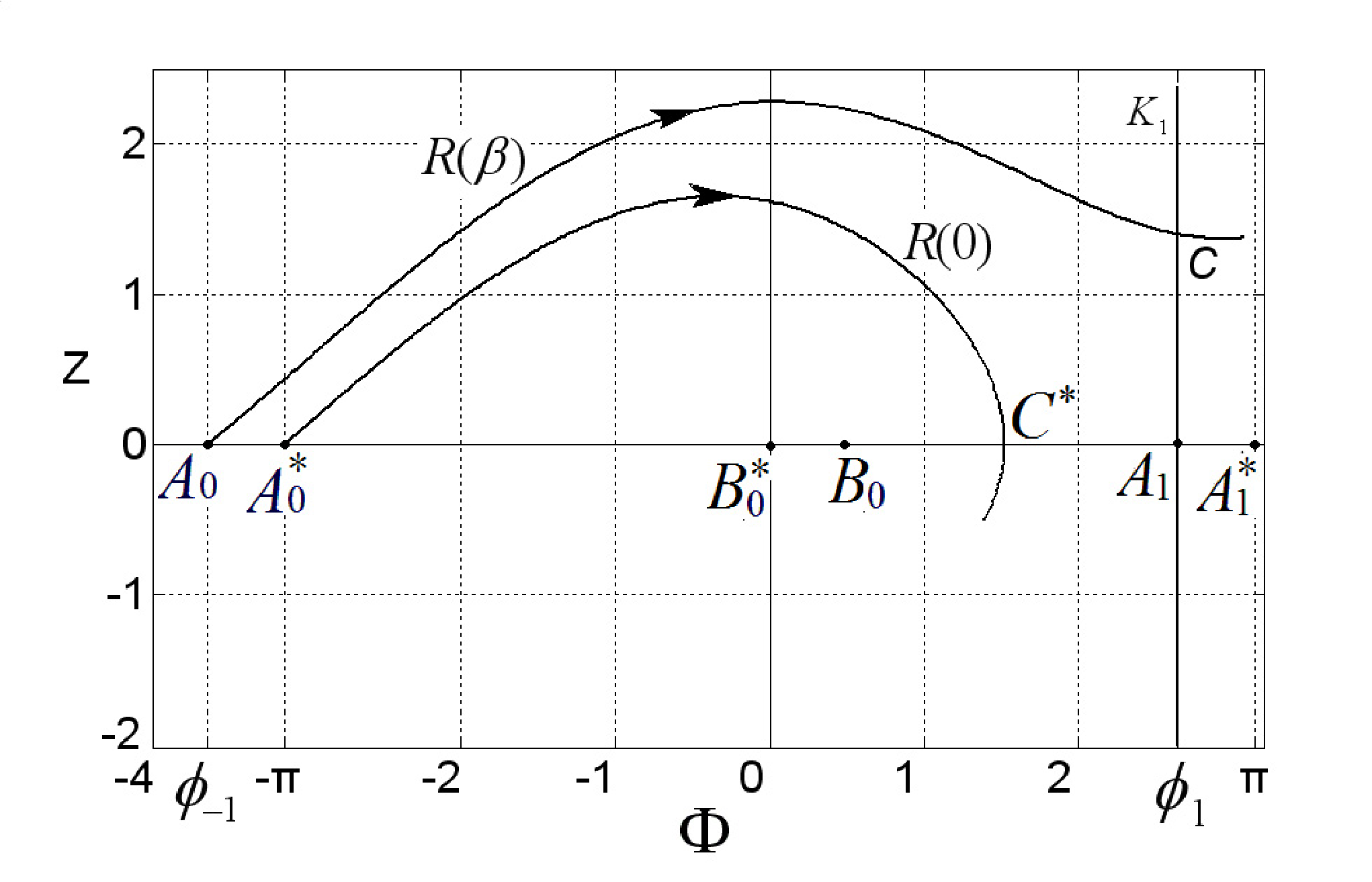}
\caption{The trajectories $R$ as $\beta=0$ and $\beta>0$.}
\label{fig:p3}
\end{figure}

\noindent
{\bf Proof.}
\par
a)\hspace*{2mm} First, Let $\beta=0$, $\Gamma =\Gamma_{1}$,
where $\Gamma_{1}>0$, thus, $\phi_{0}=0$. Three equilibria on the
main interval are denoted by $A_{0}^*$ at $(\phi_{-1}^*,0)$,
$B_{0}^*$ at $(\phi_{0}^*,0)$, $A_{1}^*$ at $(\phi_{1}^*,0)$ respectively.
According to Lemma 3, shows that the trajectory $R(0)$ emitting
from point $A_{0}^*$ must reach the point $C^*$ at $(\phi_{c}^*,0)$ on the $\Phi$-axis, and that $0\leqslant\phi_{c}^*<\pi$, as shown in Figure
\ref{fig:p3}.
\par
Now set $\beta>0$ (here $\beta\leqslant1$), thus,
$\phi_{0}=\arcsin\beta$, $0<\phi_{0}\leqslant\frac{\pi}{2}$. The
equilibria on the main interval are written as
$A_{0}$ at $(\phi_{-1},0)$, $B_{0}$ at $(\phi_{0},0)$, $A_{1}$ at $(\phi_{1},0)$
respectively. $R(\beta)$ denotes the trajectory emitting from
$A_{0}$. $K_{1}$ is a line perpendicular to the $\Phi$-axis
through the point $A_{1}$. Obviously, point $A_{0}$ is on the left
of point $A_{0}^*$, the trajectory $R(\beta)$ runs at outside of the
curve $R(0)$ at the beginning. According to Definition 2, it
follows that $\frac{\mathrm{d}\theta}{\mathrm{d}\beta}>0$ for
the vectors of field on the curve $R(0)$ except point $A_{0}^*$ and $C^*$.
This means that, while $\beta>0$, the vectors of field are rotated anticlockwise, moreover $-\frac{\pi}{2}<\theta\leqslant\frac{\pi}{2}$, therefore the direction of the vector at each
point of curve $R(0)$ is towards outside. Hence, according to the
uniqueness of solutions and Lemma 1, it follows that the
trajectory $R(\beta)$ must run at the outside of $R(0)$, and that, reaches
point $C$ at $(\phi_{c},0)$ on the $\Phi$-axis, where
$\phi_{max}\leqslant\phi_{c}\leqslant\phi_{1}$,
$\phi_{max}=max\{\phi_{c}^*,\phi_{0}\}$; or reaches point
$C$ at $(\phi_{1},z_{c})$ on the line $K_{1}$, where
$0<z_{c}<\frac{\beta+1}{\Gamma_{1}}$.
\par
Assume that $R(\beta)$ intersects point $C$ at $(\phi_{1},z_{c})$ on line
$K_{1}$, where $0<z_{c}<\frac{\beta+1}{\Gamma_{1}}$, as shown in Figure \ref{fig:p3}.
According to the continuous dependence of solutions on parameters, one
can find a trajectory $R(\beta_{0})$, where $0<\beta_{0}<\beta$,
which intersects point $A_{1}$. Proposition a) is correct.
\par
Assume that $R(\beta)$ intersects point $C$ at $(\phi_{c},0)$ on the $\Phi$-axis,
where $\phi_{max}\leqslant\phi_{c}\leqslant\phi_{1}$,
$\phi_{max}=max\{\phi_{c}^*,\phi_{0}\}$. Notice that point
$C$ at $(\phi_{c},0)$ isn't located on the left of point
$C^*$ at $(\phi_{c}^*,0)$ at least with $\beta$ increasing, point $A_{1}$
moving leftward monotonously. And while  $\beta=1$,
$\phi_{0}=\frac{\pi}{2}$, $B_{0}$ coincides with $A_{1}$, $C$ must
be between the two points. According to the these situations, it
follows that there must be a trajectory $R(\beta_{0})$, where
$0<\beta_{0}\leqslant1$, which intersects point $A_{1}$. Proposition a)
is also correct. The proof of proposition a) is complete.

\par
b)\hspace*{2mm} Proof by contradiction. Assume that
$\beta_{0}(\Gamma_{1})$ is not unique, there is
$\beta_{0}^\prime(\Gamma_{1})\neq\beta_{0}(\Gamma_{1})$, we may set
$\beta_{0}^\prime(\Gamma_{1})>\beta_{0}(\Gamma_{1})$.
\par
\begin{figure} [h]
\centering
\includegraphics[width=2.80in, height=1.79in]{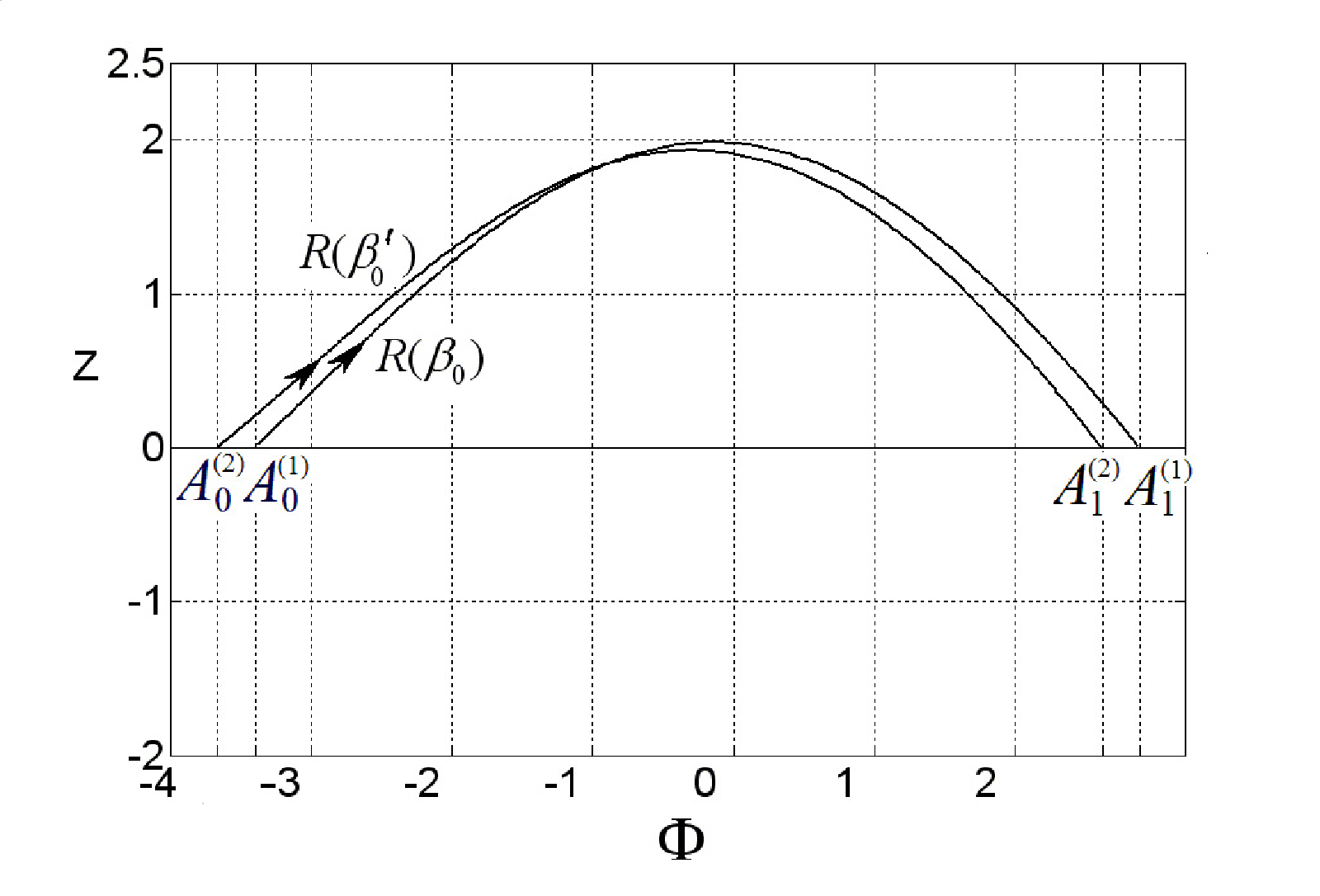}
\caption{Two critical trajectories which intersect each other.}
\label{fig:p4}
\end{figure}

While $\beta=\beta_{0}$, the two saddles are denoted by $A_{0}^{(1)}$
and $A_{1}^{(1)}$ respectively. The trajectory connecting them
is denoted by $R(\beta_{0})$, as shown in Figure \ref{fig:p4}.
\par
While $\beta=\beta_{0}^\prime$, the two saddles are denoted by
$A_{0}^{(2)}$ and $A_{1}^{(2)}$ respectively. The trajectory
connecting them is denoted by $R(\beta_{0}^\prime)$, as shown in
Figure \ref{fig:p4}.
\par
Because $\beta_{0}^\prime>\beta_{0}$, obviously, point $A_{0}^{(2)}$
is on the left of point $A_{0}^{(1)}$, point $A_{1}^{(2)}$ on the
left of point $A_{1}^{(1)}$. Therefore $R(\beta_{0}^\prime)$ must
cross with $R(\beta_{0})$, the intersection is not on the $\Phi$-axis,
as shown in Figure \ref{fig:p4}. Referring to the preceding
proofs, we know that $R(\beta_{0}^\prime)$ can not intersect
$R(\beta_{0})$ by the analysis of the rotated vector field. The
result contradicts, hence
$\beta_{0}(\Gamma_{1})$ must be unique. The Proof is complete.
\par
c)\hspace*{2mm} First we prove for the case of $\beta>\beta_{0}$.
\par
\begin{figure} [h]
\centering
\includegraphics[width=2.80in, height=1.80in]{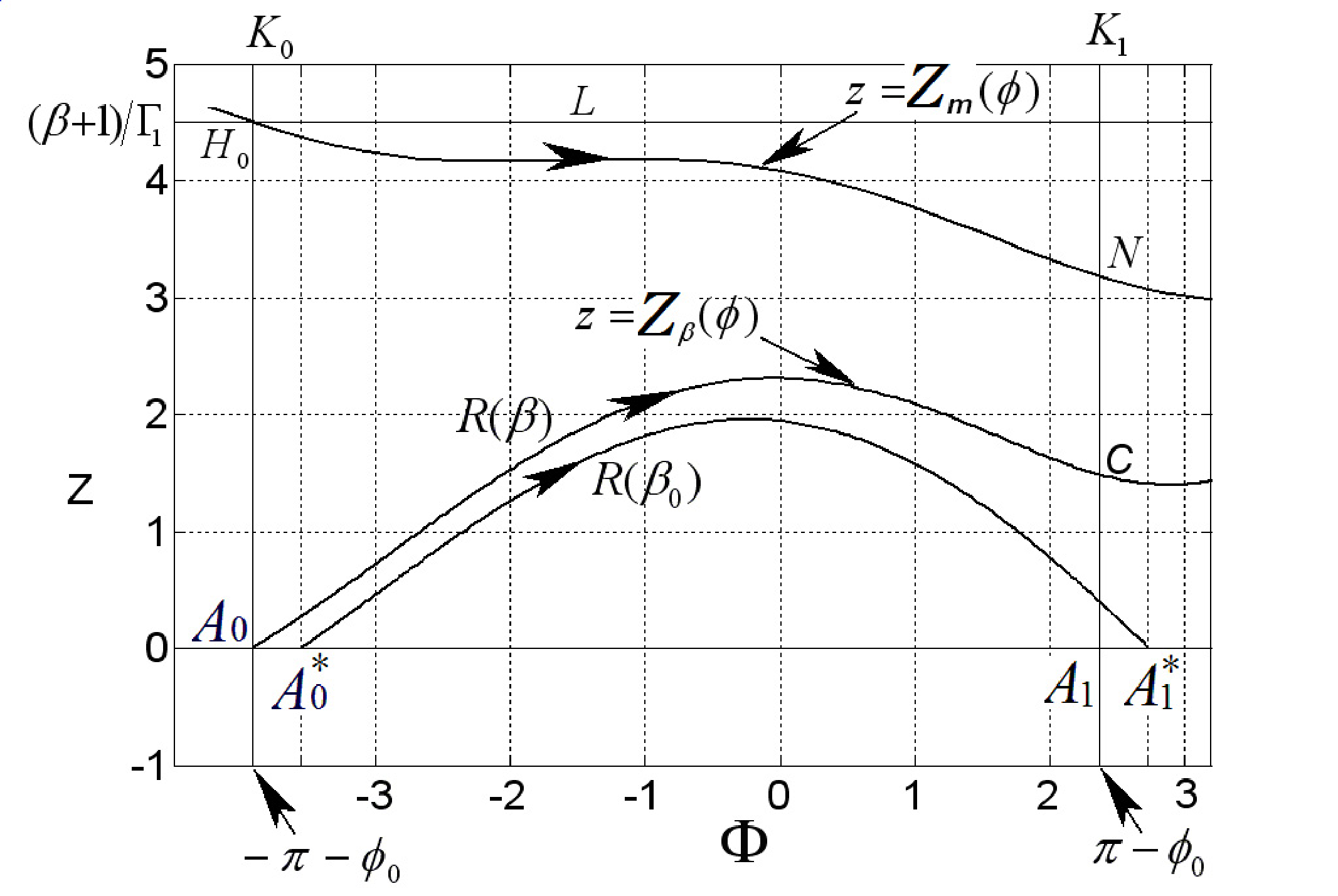}
\caption{The trajectories as $\beta>\beta_{0}$.}
 \label{fig:p5}
\end{figure}

Suppose that $\beta=\beta_{0}$, equation (\ref{eq:f2}) has a critical periodic
trajectory $R(\beta_{0})$, the saddles connecting by it are denoted by
$A_{0}^{*}$ and $A_{1}^{*}$ respectively, as shown in Figure \ref{fig:p5}.
\par
While $\beta>\beta_{0}$, the trajectory emitting from the new
equilibrium $A_{0}$ is denoted by $R(\beta)$. Referring to the
preceding proofs, it can be obtained that $R(\beta)$ must run above
$R(\beta_{0})$, as shown in Figure \ref{fig:p5}. The function for
$R(\beta)$ reads $z=Z_{\beta}(\phi)$.
Obviously, $Z_{\beta}(-\pi-\phi_{0})=0$,
$Z_{\beta}(\pi-\phi_{0})>0$, thus,
\begin{equation} \label{eq:f12}
Z_{\beta}(-\pi-\phi_{0})<Z_{\beta}(\pi-\phi_{0}).
\end{equation}
\par\noindent
Meanwhile, notice that $\displaystyle
\frac{\mathrm{d}z}{\mathrm{d}\phi}=-\frac{1+\sin\phi}{\beta+1}\cdot\Gamma_{1}\leqslant0$
for each point on the line $L$: $z=\frac{\beta+1}{\Gamma_{1}}$, and that
the vectors of field direct from the left to the right,
it follows the vectors of field on the line $L$ point downwards of
$L$. Moreover, notice that the vectors of field above the line
$L$ direct downwards absolutely. Considering this situation above, and taking a trajectory passing
through point $H_{0}$ at $(-\pi-\phi_{0},\frac{\beta+1}{\Gamma_{1}})$, for which the
function reads $z=Z_{m}(\phi)$, it follows that,
$Z_{m}(-\pi-\phi_{0})=\frac{\beta+1}{\Gamma_{1}}$,
$Z_{m}(\pi-\phi_{0})\leqslant\frac{\beta+1}{\Gamma_{1}}$, and thus,
\begin{equation} \label{eq:f13}
Z_{m}(-\pi-\phi_{0})\geqslant Z_{m}(\pi-\phi_{0}).
\hskip 2mm
\end{equation}
\par\noindent
According to inequality (\ref{eq:f12}), (\ref{eq:f13}), and the
continuous dependence of solutions on initial
conditions, on the segment $A_{0}H_{0}$, there exists a point $P$
through which the trajectory satisfies that $Z_{T}(-\pi-\phi_{0})=
Z_{T}(\pi-\phi_{0})$, obviously, $Z_{T}(\pi)>0$. Due to the
periodicity of the cylinder $E^2$, it follows that:
\begin{equation} \label{eq:f14}
Z_{T}(\phi)=Z_{T}(\phi+2\pi), \quad\quad -\infty<\phi<+\infty,
\end{equation}
\par\noindent
i.e. a periodic solution, and that
$Z_{T}(\phi)>0$, $\beta>\beta_{0}$. The proposition is correct.
\par
Second we prove for the case of $\beta<\beta_{0}$.
\par
Suppose that $\beta=\beta_{0}$, the two saddles are denoted by $A_{0}^*$
and $A_{1}^*$ respectively. Let $R(\beta_{0})$ denote the
trajectory connecting two saddles, for which the function reads
$z=Z_{0}(\phi)$, hence $Z_{0}(\phi)=Z_{0}(\phi+2\pi)$, $-\infty<\phi<+\infty$.
\par
\begin{figure} [h]
\centering
\includegraphics[width=2.80in, height=1.79in]{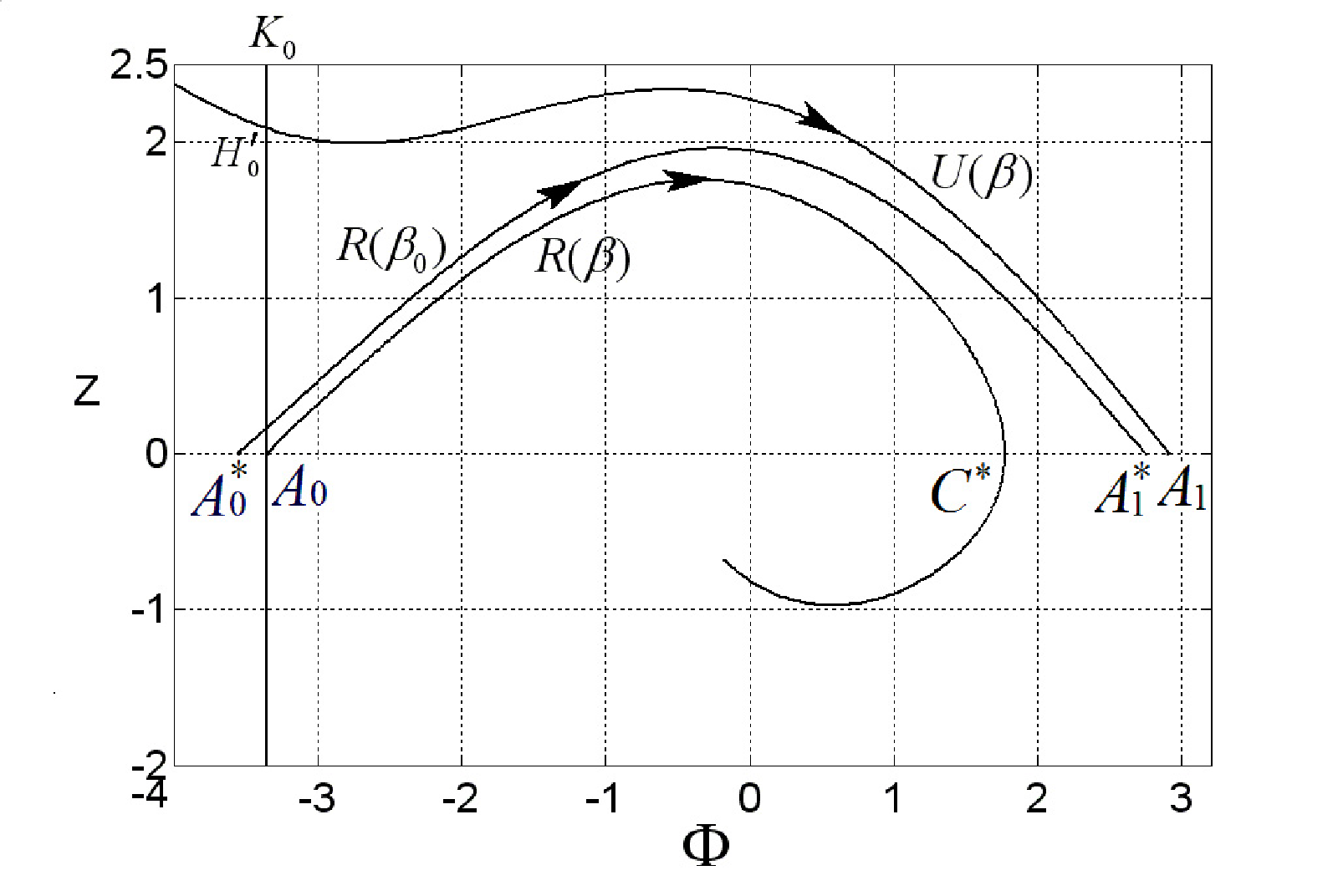}
\caption{The trajectories as $\beta<\beta_{0}$.}
 \label{fig:p6}
\end{figure}

While $\beta<\beta_{0}$, the two saddles are denoted by $A_{0}$ and
$A_{1}$ respectively. Let $R(\beta)$ denote the trajectory emitting
from point $A_{0}$, Let $U(\beta)$ denote the one coming into point
$A_{1}$, as shown in Figure \ref{fig:p6}. We will prove by
contradiction. Assuming there is a periodic solution, for which the function reads
$z=Z_T(\phi)$, obviously,
$Z_T(\phi)=Z_T(\phi+2\pi)$, $-\infty<\phi<+\infty$.
\par
Analogous to the preceding proofs, it is easy to derive the
relation of $R(\beta_{0})$ with $R(\beta)$
and $U(\beta)$ in position, i.e., $R(\beta)$, $U(\beta)$ at each
side of $R(\beta_{0})$ respectively, as in Figure \ref{fig:p6}.
Obviously, if $z=Z_T(\phi)$ passes through the
segment $A_{0}H_{0}^\prime$, it must be through the $\Phi$-axis.
Observing the directions of the vectors of field at each side of
the $\Phi$-axis, shows that the trajectory passing
through the $\Phi$-axis can not be a periodic one. Therefore,
$z=Z_T(\phi)$ is impossible to pass through the segment
$A_{0}H_{0}^\prime$. So there are only two cases of the relation
between $z=Z_T(\phi)$ and $z=Z_{0}(\phi)$ as following.
\par
In the first case, $Z_T(\phi_{-1})<Z_{0}(\phi_{-1})$, because
the trajectory $z=Z_T(\phi)$ doesn't crossover the
$\Phi$-axis, we have $Z_T(\phi)\leqslant 0$,
$-\infty<\phi<+\infty$. While $\beta>0$, $\Gamma>0$, according
to Lemma 2 shows that this kind of periodic solution is
impossible.
\par
In the second case, trajectory $z=Z_T(\phi)$ is above $R(\beta_0)$, i.e. $Z_T(\phi)>Z_{0}(\phi)$,
$\phi_{-1}^*\leqslant\phi\leqslant\phi_{1}^*$. While $\beta<\beta_{0}$,
according to Lemma 2, shows that this situation is also
impossible.
\par
Summing up the above conclusions gives that equation (\ref{eq:f2})
does not exists the periodic solution while $\beta<\beta_{0}$. The
proof is complete for the proposition of $\beta<\beta_{0}$.
\begin{figure} [h]
\centering
\includegraphics[width=1.40in, height=1.60in]{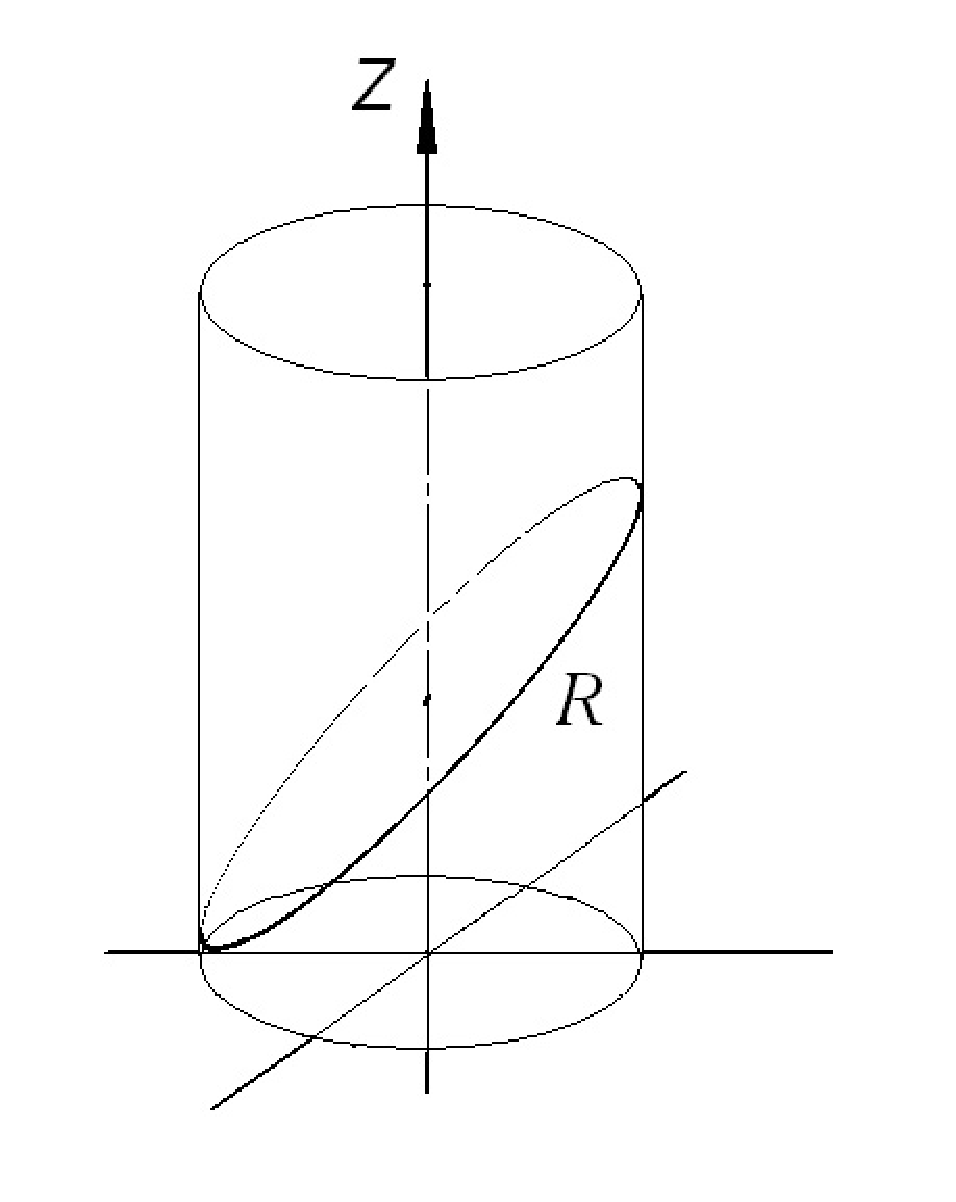}
\caption{The limit cycle on the cylinder $E^2$.}
 \label{fig:p7}
\end{figure}
\par
d)\hspace*{2mm} According to Lemma 2, the periodic solution is unique obviously,
and is the limit cycle of type II on $E^2$, as shown in Figure
\ref{fig:p7}. The average integral of the divergence of flow over the cycle satisfies
\begin{equation}\label{eq:f15}
\frac{1}{T}\int_0^T\triangledown\cdot\textbf{F}(\phi,z) \mathrm{dt}=-\Gamma<0,
\hskip 3mm
\end{equation}
therefore it is a stable limit cycle \cite{Perk}. The periodic
solution is a stable periodic solution.
\par
e)\hspace*{2mm} Suppose $\beta_{0}(\Gamma_{1})=1$,
while $\beta=\beta_{0}(\Gamma_{1})$, $\Gamma=\Gamma_{1}$, there is a critical
periodic trajectory connecting two equilibria $A_{0}$ and $A_{1}$, denoted by $R(\beta,\Gamma_{1})$,
as shown in Figure \ref{fig:p8}. Referring to the preceding proof, consider $\Gamma_{2}>\Gamma_{1}$,
we conclude that, while $\Gamma=\Gamma_{2}$ and $\beta$ unchanged, the
trajectory emitting from $A_{0}$, denoted by $R(\beta,\Gamma_{2})$, must run at the
inside of $R(\beta,\Gamma_{1})$, as shown in Figure \ref{fig:p8}. According to Lemma 1,
and noticing that the equilibrium $B_{0}$ coincides
with $A_{1}$ as $\beta=1$, it follows that $R(\beta,\Gamma_{2})$
must intersect the $\Phi$-axis at point $A_{1}(i.e.$ $B_{0})$,
it is a critical periodic trajectory yet. Notice that
$\beta=1$ still. According to the uniqueness of $\beta_{0}$, it
follows that $\beta_{0}(\Gamma_{2})=1$, while $\Gamma_{2}>\Gamma_{1}$. In
this case, the slope by which $R(\beta,\Gamma_{2})$ approaching the point
$A_{1}$ is $\lambda_{1}=0$ .
    All proofs for this theorem are complete.
\begin{figure} [h]
\centering
\includegraphics[width=2.80in, height=1.79in]{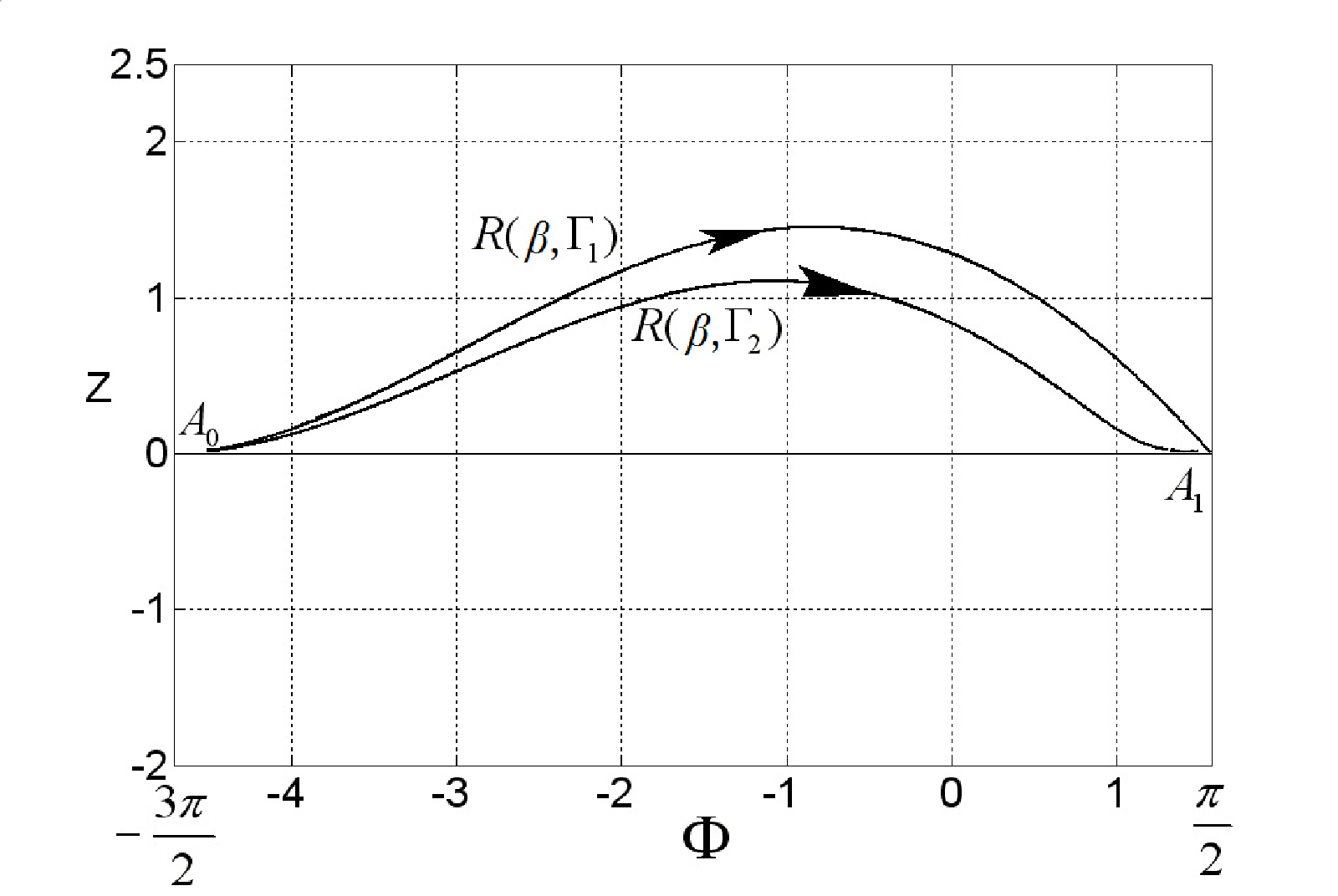}
\caption{The critical periodic trajectory as
$\Gamma_{2}>\Gamma_{1}$.}
 \label{fig:p8}
\end{figure}

\begin{flushleft}
\large \bf 4.\hspace*{2mm} Discussion of the theorem
\end{flushleft}

According to the theorem above, it follows a conclusion that, for any
$\Gamma>0$, equation (\ref{eq:f2}) exists a critical value
$\beta_{0}$. In other words, while $\beta\geqslant\beta_{0}$, it exists a
periodic solution that is unique, stable and positive. It indicates
that, whatever the initial state it is, the final state will be in a
self-sustained periodic motion toward positive. While $\beta<\beta_{0}$,
equation (\ref{eq:f2}) exists no any periodic solution. In this situation, all
trajectories on the phase plane will come into the saddles or
focuses ultimately. The focuses are stable but saddles
unstable. If there are some disturbances, it will deviate from the
saddles, eventually, will stay at the
focuses for ever. In physics, whatever the initial state it is, the system
will be static at the focus after undergoing sufficient time.
\par
According to Proposition e) it has been shown that, if there
exists a minimum value $\Gamma_{min}$ such that $\beta_{0}(\Gamma_{min})=1$, while
$\Gamma\geqslant\Gamma_{min}$, the critical value $\beta_{0}$ will be 1
all the time. By the numerical calculations, M.Urabe has given the
definite value \cite{Urabe}, that is $\Gamma_{min}\simeq1.193$. From
the experimental values given by Portis in CDW's materials \cite{Portis}, it follows that
$\Gamma=5.83$. Clearly, the
practical value of $\Gamma$ is large than the minimum value
$\Gamma_{min}$ so far. So we can believe that, in the practical
conditions, $\beta_{0}\equiv1$.
\par
These conclusions above are employed in the research of CDW. The
derivations in detail was presented in other paper \cite{Li}. It is
unnecessary to go into details in this paper.
\par
\vspace{0.4cm} \noindent {\large \bf {Acknowledgments}}
\par The author thanks the referee for his/her
constructive and useful recommendations.


\end{document}